\documentclass{article}

\usepackage[english]{babel}
\usepackage[utf8]{inputenc}
\usepackage{amsmath,amssymb}
\usepackage{parskip}
\usepackage{graphicx}
\usepackage[style=alphabetic]{biblatex}
\usepackage[top=2.5cm, left=3cm, right=3cm, bottom=4.0cm]{geometry}
\usepackage[table]{xcolor}

\usepackage{amssymb}
\usepackage[T2A]{fontenc}  
\usepackage[utf8]{inputenc} 
\usepackage[english]{babel} 
\usepackage{mathtools}

\usepackage{hyperref}

\usepackage{stmaryrd}

\usepackage{amsmath, amssymb, amsthm}

\providecommand{\openbox}{\leavevmode
  \hbox to.77778em{%
  \hfil\vrule
  \vbox to.675em{\hrule width.6em\vfil\hrule}%
  \vrule\hfil}}

\providecommand{\keywords}[1]{\textbf{Keywords:} #1}

\title{Random walks on rank one symmetric spaces \\ of noncompact type$^*$}
\author{Fedor Gnetov$^\dagger$, Valentin Konakov$^\ddagger$}
\date{}

\begin{document}
\maketitle

\section*{Abstract}

We establish a central limit theorem, a local limit theorem, and a law of large numbers for a natural random walk on a symmetric space \(M\) of non-compact type and rank one. This class of spaces, which includes the complex and quaternionic hyperbolic spaces and the Cayley hyperbolic plane, generalizes the real hyperbolic space \(\mathbb{H}^{n}\). Our approach introduces a unified algebraic framework that generalizes the Möbius addition, previously used for the constant curvature case, to define the random walk via a non-Euclidean summation of variables. We demonstrate that the renormalized walk converges to the heat kernel associated with the Laplace-Beltrami operator on \(M\), which plays the role of the limiting normal law. The proofs leverage the harmonic analysis of spherical functions on symmetric spaces. To the best of our knowledge, these results are new in the context of rank one symmetric spaces.

\keywords{central limit theorem, local limit theorem, law of large numbers, random walk, symmetric spaces of rank one, hyperbolic spaces, heat kernel, Laplace-Beltrami operator, harmonic analysis, spherical functions, M\"obius addition.}

\section{Introduction}

Let \(M\) denote a symmetric space of non-compact type and rank one. This class of spaces, which includes the complex and quaternionic hyperbolic spaces and the Cayley hyperbolic plane, generalizes the real hyperbolic space \(\mathbb{H}^{n}\) (\(n\geqslant 2\)). These are complete, simply connected Riemannian manifolds with strictly negative, yet non-constant, sectional curvature.

Our objective is to establish limit theorems—specifically, a central limit theorem, a local limit theorem, and a law of large numbers—for a natural random walk on \(M\). In preceding work \cite{KonMen25}, we addressed this problem for the Poincaré hyperbolic space (the rank one space with constant curvature). Therein, building upon the foundational work of Karpelevich et al. \cite{Karpelevich59}, we showed that a suitably renormalized sum of random variables, defined via the Möbius addition on the ball model, converges in law to the fundamental solution of the associated heat equation. The limit was identified as the normal law on \(\mathbb{H}^{n}\), which corresponds to the heat kernel defined by the Laplace-Beltrami operator. 

The approach in \cite{Karpelevich59} and \cite{KonMen25} relied on the Möbius addition to provide a non-Euclidean analogue for summing random variables. In more detail, the connection with harmonic analysis on underlying geometric structures provides a powerful framework for defining random walks and  \\
\rule{5cm}{1pt} 

*The study has been funded by the Russian Science Foundation (project № 24-11-00123). \\
\begin{tabular}{@{}l}
$\dagger$National Research University Higher School of Economics, Russian Federation / HSE University, RF \\
$\ddagger$National Research University Higher School of Economics, Russian Federation / HSE University, RF
\end{tabular}

investigating their limit theorems, as it offers natural generalizations of characteristic functions and variances; see \cite{Terras2013} and \cite{helgason1984groups} for results in this direction. A more functional perspective is provided by the works of Pinsky \cite{Pinsky75}, \cite{Pinsky76}, who established the convergence of the isotropic transport process to Brownian motion on manifolds with Ricci curvature bounded below. For a discrete space-time approximation in the Poincaré plane, we refer to Gruet \cite{Gruet08}.

In the present work, we generalize this program to all rank one symmetric spaces. A key innovation is the introduction of a unified algebraic framework that generalizes the Möbius addition, allowing for a coherent definition of a random walk on \(M\) via a non-Euclidean "sum" of variables; this framework is rooted in the classical theory of symmetric spaces, for which we primarily follow Helgason's treatises \cite{helgason1979differential}. The connection between such algebraic operations, probability, and harmonic analysis is further elaborated in the monograph of Terras \cite{Terras2013}. For the specific harmonic analysis required to analyze this walk—providing the natural generalizations of Fourier transforms and characteristic functions—we rely fundamentally on Helgason's presentation \cite{helgason1984groups} of the theory of spherical functions and related constructs, which subsumes the approach used by Platonov \cite{Platonov99}.

Within this general setting, we demonstrate that the renormalized random walk converges to the heat kernel, just as in the specific hyperbolic case. This heat kernel, defined as the fundamental solution of the heat equation involving the Laplace-Beltrami operator on \(M\), plays the role of the limiting normal law. For a detailed analysis of the heat kernel on symmetric spaces and its properties, we refer to the comprehensive work of Anker \cite{anker:hal-00002509}. 

Whereas local limit theorems have been previously established for some nilpotent groups, such as the Heisenberg group (see Breuillard \cite{Breuillard05}), our results in the present context of rank one symmetric spaces are, to the best of our knowledge, new. Our proof of the local limit theorem adapts and extends the strategy we developed in \cite{KonMen25} for the constant curvature case, leveraging the deeper algebraic structure now available. As a direct by-product of our analysis, we also obtain a law of large numbers.

Our work can be compared with that of Kogler \cite{Kogler25}, who also studied limit theorems on symmetric spaces. The advantages of our work are: 1. the correspondence between geometry and algebraic structure, since we introduce a generalization of Möbius addition; 2. we establish different types of convergence.

The paper is organized as follows. We describe the generalized algebraic structure, the corresponding random walk, and our main convergence results in Section 2. Section 3 is devoted to the necessary tools from harmonic analysis on symmetric spaces. Finally, the proofs of the main theorems are presented in Section 4.

\section{Main results}

\subsection{Setting}

Let \( M \cong G/K \) be a symmetric space of noncompact type, where \( G \) is a connected semisimple Lie group with finite center. Fix an arbitrary point $p \in M$ on the manifold $M$. \( K \subset G \) is a maximal compact subgroup. Let \( \theta \) be the corresponding Cartan involution, with \( \mathfrak{g} = \mathfrak{k} \oplus \mathfrak{p} \) the Cartan decomposition of the Lie algebra, where $\mathfrak{k}$ and $\mathfrak{p}$ are the $+1$ and $-1$ eigenspaces for $\theta$. 

Fix the Iwasawa decomposition  
\[
G = KAN, \quad \mathfrak{g} = \mathfrak{k} \oplus \mathfrak{a} \oplus \mathfrak{n},
\]  
where \( \mathfrak{a} \subset \mathfrak{p} \) is a maximal abelian subspace and \( \mathfrak{n} \) is the nilpotent subalgebra associated with a choice of positive roots. For \( g \in G \), define \( H(g) \in \mathfrak{a} \) and \( n(g) \in N \) uniquely by the relation  
\[
g = \kappa(g) \exp(H(g)) \, n(g), \quad \kappa(g) \in K.
\]  
Let \( \Sigma \) be the set of roots of \( \mathfrak{g} \) with respect to \( \mathfrak{a} \), \( \Sigma^+ \) the set of positive roots corresponding to the choice of \( \mathfrak{n} \), and let \( m_\alpha \) denote the multiplicity of a root \( \alpha \). The half-sum of positive roots, counted with multiplicity, is defined as  
\[
\rho = \frac{1}{2} \sum_{\alpha \in \Sigma^+} m_\alpha \alpha.
\]  
Since we observe a rank 1 symmetric space, $\Sigma^+ = \{\alpha, 2 \alpha \}$. Take $H_0 \in \mathfrak{a}$ to satisfy condition $\alpha(H_0) = 1$.  Let \( \mathfrak{a}_{\mathbb{C}}^* \) be the space of complex-linear forms on \( \mathfrak{a} \).
Let \( Z \) and \( N \) denote the centralizer and normalizer of \( \mathfrak{a} \) in \( K \), respectively. The quotient \( N/Z \) is the Weyl group \( W = W(\mathfrak{g}, \theta) \). Let \( A = \exp \mathfrak{a} \), and let \( A^+ = \exp \mathfrak{a}^+ \), where \( \mathfrak{a}^+ \) is a positive Weyl chamber. The real rank of \( \mathfrak{g} \) and \( G \) is defined as \( d = \dim \mathfrak{a} \). Fix $KA^+K$ decomposition:

$$g = k_1 \exp(c(g)) k_2, \quad k_1, k_2 \in K, \ c(g) \in a^+$$

The Killing form \( B: \mathfrak{g} \times \mathfrak{g} \to \mathbb{R} \) is a symmetric, non-degenerate bilinear form defined by  
\[
B(X, Y) = \operatorname{tr}(\operatorname{ad}(X) \circ \operatorname{ad}(Y)),
\]  

We define a positive definite inner product on $\mathfrak{g}$ by

\[
\langle U,V \rangle := -B(U, \theta V), \ \| U \| = \sqrt{\langle U,U \rangle}    \quad (U,V \in \mathfrak{g}).
\]


Let $dg$ and $dk$ denote the Haar measures on $G$ and $K$, respectively, with $dk$ normalized such that $\int_K dk = 1$. Let $\mu_M$ denote the Riemannian volume measure on $M$.

Let $\pi$ denote the natural mapping $g \mapsto g \cdot p$ of $G$ onto $M$.

Let \( x \) and \( y \) be points of \( M \) and we want to define \( x\oplus y \in M \) the sum of these points. The construction of the sum will consist of three stages:


1) Using (Theorem 3.3(i) \cite{helgason1979differential}) and (Theorem 3.3(iii) \cite{helgason1979differential}) we can construct geodesics with conditions:
\( \gamma_{(d\pi)_e X}(t) = \exp(tX) \cdot p \), \( x = \exp(X) \cdot p \), 
\( \gamma_{(d\pi)_e Y}(t) = \exp(tY) \cdot p \), \( y = \exp(Y) \cdot p \), \\
where \( X \in \mathfrak{p} \) and \( Y \in \mathfrak{p} \).

2) We make a parallel translate of the vector \((d\pi)_e Y \in T_p M\) along the geodesic \( \gamma_{(d\pi)_e X}(t) \) using for each \( t \geq 0 \) the differential \((d \exp tX)_{p} ((d\pi)_e Y)\). Consider \( Y' := (d \exp X)_{p} ((d\pi)_e Y) \in T_x M \).

3) We release the geodesic from the point \( x \) with a velocity vector \( Y' \) and as the sum \( x\oplus y \), we take the point at which the geodesic arrives at moment 1.

\textbf{Lemma 1.1.} \textit{The geodesic referred to in point (3) of the definition of addition has the form}
\[
\tilde{\gamma}(t) = \exp(X) \circ \exp(tY) \cdot p
\]
\textit{and, hence, \( x\oplus y = \exp(X) \circ \exp(Y) \cdot p \).}

\begin{proof}
We have to show that \(\tilde{\gamma}(t)\) is geodesic starting from \( x = \exp(X) \cdot p \) with initial velocity \( Y' \). Note that \(\exp(X) \in G\), that is isometry, and \(\exp(tY) \cdot p\) is geodesic according to (Theorem 3.3(iii) \cite{helgason1979differential}). So \(\exp(X) \cdot \exp(tY) \cdot p\) is geodesic as an image of a geodesic in isometry. Next,
\[
\tilde{\gamma}(0) = \exp(X) \circ \exp(0 \cdot Y) \cdot p = \exp(X) \cdot p_0 = x,
\]
\[
\frac{d\tilde{\gamma}(t)}{dt}\bigg|_{t=0} = \frac{d (\exp(X) \circ \exp(tY) \cdot p)}{dt}\bigg|_{t=0} = 
\]
\[
d(\exp(X))_p \left[\frac{\exp(tY) \cdot (p)}{dt}\right]\bigg|_{t=0} = d(\exp(X))_p[(d\pi)_e Y] = Y'.
\]
  
\end{proof} 

Let introduce multiplication by constant element in symmetric space. $g \in G, \ g = k_1 \exp(c(g)) k_2$.
$$v \otimes g = k_1 \exp(vc(g)) k_2, \quad v\in \mathbb{R}$$

\subsection{The random walk and associated convergence results}

Let \(Z\) be a \(M\)-valued random variable defined on some probability space \((\Omega,\mathcal{A},\mathbb{P})\). We will further assume that:

\textbf{[R]} \(Z\) has \emph{radial} density \(f_{Z}\in C^{\infty}_{0}(K \backslash G / K)\) w.r.t. the Riemannian volume.

Let now \((Z^{j})_{j\geqslant 1}\) be a sequence of i.i.d random variables which have the same law as \(Z\). Define then for \(N\in\mathbb{N}\),

\begin{equation*}
\bar{S}_{N} :=\bigoplus_{j=1}^{N}\frac{1}{N}\otimes Z^{j},
\end{equation*}
\begin{equation*}
S_{N} :=\bigoplus_{j=1}^{N}\frac{1}{\sqrt{N}}\otimes Z^{j}. 
\end{equation*}

Then the following results hold:

\textbf{Theorem 2.1 (Law of large numbers).} Under [R],
\[
\bar{S}_{N}\xrightarrow{\mathbb{P}}0.
\]

It will be shown in Section 4 that \(S_{N}\) has a density \(f_{S_{N}}\), which can be expressed as the non Euclidean convolution of the densities of the \((Z^{j})_{j\in\llbracket 1,N \rrbracket}\). We quantify below the asymptotic behavior of that density.

\textbf{Theorem 2.2 (Central limit theorem).} Under [R], it holds that for measurable sets \(A\) in $M$,
\[
\int_{M} \mathbb{I}_{A}(x)f_{S_{N}}(x)\mu_{M}(dx) \xrightarrow[N]{} \int_{M} \mathbb{I}_{A}(x)\Psi(t,x)\mu_{M}(dx),
\]
with \(t:=\frac{1}{n \| \alpha \|^2}\int_{0}^{+\infty}\tilde{\eta}^{2}\mu_{Z,R}(d\tilde{\eta})\), where \(\mu_{Z,R}\) stands for the measure induced by the law of \(Z\) in geodesic coordinates.

Also, \(\Psi(t,\cdot)\) stands for the hyperbolic heat kernel in M. Namely it denotes the fundamental solution of the equation
\[
\frac{1}{2}\Delta_{M}\Psi(t,x)=\partial_{t}\Psi(t,x),\ \Psi(0,\cdot)=\delta(\cdot).
\]

The specific expression of \(\Psi(t,\cdot)\) will be given in Section 3.2 below. It plays the same role in the current setting as the normal density in the classical Euclidean central limit theorem.

We are furthermore able to specify the previous result giving a convergence rate.

\textbf{Theorem 2.3 (Local limit theorem).} Under [R], there exists \(C:=C(m_{\alpha}, m_{2\alpha}, \mu_{Z})\) s.t. for all \(x\in M\) and \(N\) large enough
\begin{equation}
|f_{S_{N}}(x)-\Psi(t,x)|\leqslant\frac{C}{N} \Big{(}\frac{1}{t^{2}}\Big{(}\frac{1}{ t^{\frac{1}{2}}}\wedge\frac{1}{t^{\frac{n}{2}}}\Big{)}+1\Big{)}. \tag{2.1}
\end{equation}

From now on we will denote by \(C\) a generic constant, that may change from line to line and depends on $m_{\alpha}, m_{2\alpha}, t$ and the law \(\mu_{Z}\), i.e. \(C=C(m_{\alpha}, m_{2\alpha}, \mu_{Z}) \).

\section{Symmetric space: spherical transform and heat kernel}

\subsection{Spherical transform and its inverse}

For radial function $f \in C_{0}(K \backslash G / K)$ bi-K-invariant, $\lambda \in \mathfrak{a}^*$.

It's spherical transform is given by the expression (For details see \cite{helgason1984groups} page 449).
\[
\tilde{f}(\lambda) = \int_G f(g)\phi_{-\lambda}(g)  dg = C \int_{0}^{\infty} f(\exp (tH_0)) \phi_{-\lambda}(\exp (tH_0)) \sinh(t) \sinh (2t)  dt \tag{3.1}
\]

where the radial functions $\varphi_{\lambda}$ are eigenfunctions of the Laplace--Beltrami operator in $G$ expressed in radial coordinates. Namely, for a generic smooth enough $\phi: G \to \mathbb{R}$,

\[
\Delta_M \phi = \frac{\partial^2}{\partial r^2} \phi + \frac{1}{A(r)} \frac{dA}{dr} \frac{\partial}{\partial r} \phi + L_S \phi,
\]  
where \( L_S \) is the Laplace-Beltrami operator on \( S_r(p) \). Here \( 0 < r < \infty \).

$$A(r) = C \sinh^{m_{\alpha}}(cr) \sinh^{m_{2 \alpha}} (2cr)$$

where $c = (2m_{\alpha} + 8m_{2 \alpha})^{-1/2}$

When $\phi$ is radial, this operator rewrites as

$$\Delta_M \phi = \frac{\partial^2}{\partial r^2} \phi + \frac{1}{A(r)} \frac{dA}{dr} \frac{\partial}{\partial r} \phi $$

because radial functions dont depends of angle. So $L_S$ is equal to 0 (for details see \cite{helgason1984groups} page 313). And the function \(\varphi_\lambda\) in (3.1) solves the differential equation:

\begin{equation}
\begin{cases} 
\Delta_M \varphi_\lambda(r) + (\| \lambda \|^2 + \| \rho \|^2)\varphi_\lambda(r) = 0, \\ 
\varphi_\lambda(0) = 1. \tag{3.2}
\end{cases}
\end{equation}




We then define the corresponding elementary spherical function \(\varphi_\lambda\) setting for all \(g \in G, \ \lambda \in \mathfrak{a}^* \),

\[\varphi_\lambda(g) = \int_K e^{(i\lambda - \rho)(H(gk))}  dk, \tag{3.3}\]

\textbf{Lemma 3.1.}
$$\frac{\partial^2}{\partial t^2} \varphi_{\lambda} \left(\exp \left(\frac{t H_0}{\| H_0 \|} \right) \right) \biggr|_{t=0} =  -\frac{\| \lambda\|^2 + \| \rho \|^2}{n}, \quad t \in \mathbb{R}_+$$
\begin{proof}
    Obtained by direct calculation by substituting the Taylor expansion of $\varphi_\lambda$ into (3.2).
\end{proof}

\subsection{The heat kernel on $M$}

The heat kernel will provide the limit law which is somehow the analogue in the current non-Euclidean setting of the normal law. The normal density of parameter \(t>0\) in \(M\) is defined as the solution to

\[
\frac{1}{2}\Delta_{M}\Psi(t,x)=\partial_{t}\Psi(t,x),\ \Psi(0,\cdot )=\delta(\cdot).
\]

In the literature the usual heat equation considered is (see, e.g., \cite{anker:hal-00002509})

\[
\Delta_{M}\psi(t,x)=\partial_{t}\psi(t,x),\ \psi(0,\cdot)=\delta(\cdot).
\]

It can actually be solved through the spherical transform (assuming the solution is radial). Namely, from (3.1), one derives that for all \(\lambda\in \mathfrak{a}^*\):

\begin{equation}
\hat{\psi}(t,\lambda)=\exp(-(\| \lambda \|^{2}+ \| \rho \|^{2})t). \tag{3.4}
\end{equation}




From the previous definitions it is clear that for $\Psi$ as in Theorem 2.2 $\Psi(t, x) = \psi(\frac{t}{2}, x)$.

For $\mathbb{H}^n$ heat kernel has an implicit formula (in polar coordinates) for $m \in \mathbb{N}$:

\[
\psi(t,x) = \psi\left( t, \text{atanh} \left( \frac{\eta}{2} \right) \right)
= 
\begin{cases}
\frac{\exp(-m^2t/2)}{(2\pi)^m\sqrt{2\pi t}} \left( -\frac{1}{\sinh\eta}\partial_\eta \right)^m \exp\left(-\frac{\eta^2}{2t}\right), & n = 2m + 1, \\
\frac{\exp(-(m-1/2)^2t/2)}{(2\pi)^m\sqrt{\pi t}} \int_{\eta}^{+\infty} \frac{ds}{\sqrt{\cosh(s) - \cosh(\eta)} } \\
\quad \times (-\partial_s) \left( -\frac{1}{\sinh s}\partial_s \right)^{m-1} \exp\left(-\frac{s^2}{2t}\right), & n = 2m.
\end{cases}
\]

\subsection{Some additional tools from harmonic analysis on M.}

\subsubsection{Convolution}

\textbf{Theorem 3.1.} 
Since \( K \) is compact, there exists a unique \( G \)-invariant Radon measure \( \mu \) on \( X \) (induced by the Riemannian volume form) such that the following integral formula holds (see, e.g., \cite{knapp2013lie} page 534, Theorem 8.36, page 538; \cite{helgason1984groups} page 90-91):

\[
\int_G f(g)  dg = \int_{G/K} \left( \int_K f(gk)  dk \right) d\mu(gK), \quad \text{for } f \in C_0(G).
\]

Let \( \phi, \psi \in C_0(G) \) be \( K \)-bi-invariant functions. They may be identified with functions on \( M \) via \( \phi(gK) = \phi(g) \). We define two notions of convolution:

- The convolution on \( M \) is defined by:

$$(\phi *^M \psi) (y) = \int_{X} \phi(-x \oplus y) \psi(x) \mu(dx), \quad y \in X,$$

where $\mu$ is the Riemannian volume (measure).

- The convolution on \( G \) is defined by:

$$(\phi *^G \psi) (h) = \int_{G} \phi(g^{-1}h) \psi(g) dg, \quad h \in G. $$

The following theorem relates these two convolutions.

\textbf{Theorem 3.2.} Let \( \phi, \psi \in C_0(K \backslash G / K) \) be \( K \)-bi-invariant. For \( y = g_y \cdot p \in M \), we have
$$(\phi *^M \psi) (y) = \operatorname{vol}(K, dk) (\phi *^G \psi) (g_y)$$

\begin{proof}

    Consider the function \( F: G \times G \to \mathbb{R} \) defined by

    $$F(g, h) = \phi(g^{-1} h) \psi(g)$$
    
    We first show that \( F \) is right \( K \)-invariant in both arguments. For \( k \in K \):
    
    $$F(gk, h) = \phi((gk)^{-1} hk) \psi(gk) = \phi(k^{-1}g^{-1} h) \psi(gk) = \phi(g^{-1} h) \psi(g) =  F(g, h),$$
    $$F(g, hk) = \phi(g^{-1} hk) \psi(g) = \phi(g^{-1} h) \psi(g)= F(g, h),$$
    where we used the \( K \)-bi-invariance of \( \phi \) and \( \psi \). Thus, \( F \) descends to a well-defined function on \( M \times M \).

    Let \( y = g_y \cdot p \). And since we can select any representative of $gK, g_yK$, let's take $\exp(H_g), \  \exp(H_y)$, where $H_g, \ H_y \in \mathfrak{p}$, $\exp(H_g)K = gK, \ \exp(H_y)K = g_yK$, $H_g,  H_y$ are unique (see (Theorem 3.3(iii) \cite{helgason1979differential})).

    Using the integral decomposition (Theorem 3.1) and the invariance of \( F \):

    \[
    \int_{G} \phi(g^{-1}g_y) \psi(g) dg = \int_G F(g, g_y)  dg = \int_{G/K} \left[ \int_K F(gk, g_y)  dk \right] d\mu(gK)
    \]
    \[
    = \int_K  dk \int_{G/K} F(\exp(H_g), \exp(H_y))  d\mu(gK)
    \]
    \[
    = \int_K  dk \int_{G/K} \phi(\exp(H_g)^{-1} \exp(H_y)) \psi(\exp(H_g))  d\mu(gK) 
    = C \int_{M} \phi(-x \oplus y) \psi(x) \mu(dx)
    \]
\end{proof}

\subsubsection{Some useful inequalities}

\textbf{Lemma 3.2:}
$$\| H(\exp(t H_0)k) \| \leq  t\|H_0\|, \quad k\in K, \ t \in \mathbb{R}_+$$
$$| \alpha(H(\exp(t H_0)k)) | \leq  t\|H_0\| \|H_{\alpha}\|, \quad k\in K, \ t \in \mathbb{R}_+$$

\begin{proof}
The main result of \cite{kostant1973}, known as `Kostant's (nonlinear) convexity theorem' characterizes the image under $H$ of the set $aK$, for $a \in A$, as follows (see, e.g. \cite{balibanu2015convexitytheoremssemisimplesymmetric}):

\[
H(aK) = \text{conv}(W(\mathfrak{a}) \cdot \log a).
\]

Here `conv' indicates that the convex hull in $\mathfrak{a}$ is taken. Let's write it for $H(\exp(t H_0)k)$ in analytic way

$$H(\exp(t H_0)k) = \sum_i c_i Ad(k_i)(tH_0),$$

where $\sum_i c_i = 1$ and $k_i \in K$ correspond to the i-th element in Weyl group. Taking the norm of both sides and applying the triangle inequality yields:

$$\| H(\exp(t H_0)k) \| = \| \sum_i c_i Ad(k_i)(tH_0) \| \leq  \sum_i c_i \| Ad(k_i)(tH_0) \|$$

From Helgason \cite{helgason1979differential} page 253 Lemma 1.2: $Ad(k)$ is isometry, thus we can simplify inequality:

$$\| H(\exp(t H_0)k) \| \leq  \sum_i c_i \| tH_0 \| = \| tH_0 \|\sum_i c_i = t\|H_0\|.$$

Using Cauchy-Schwarz and the previous inequality:

$$| \alpha(H(\exp(t H_0)k)) | = \left| \langle H_\alpha, H(\exp(t H_0)k) \rangle \right| \leq t\|H_0\| \|H_{\alpha}\|.$$
    
\end{proof}

\textbf{Lemma 3.3:} Since we observe rank 1 space, we can parametrize $\mathfrak{a}^*$ using $s \in \mathbb{R}$ as $s \alpha$, $t \in \mathbb{R}_+$.

\[ \left| \frac{\partial^4}{\partial s^4} \varphi_{s}(\exp(t H_0))\right| \leq C t^2 \left| \frac{\partial^2}{\partial s^2} \varphi_{s} (\exp(t  H_0)) |_{s=0} \right|, \tag{3.5} \]

\[
-\frac{\partial^2}{\partial \lambda^2} \phi_{\lambda \alpha}(\exp(tH_0)) \leq (t \| H_0 \|  \| H_\alpha \|)^2 \varphi_{0} (\exp(tH_0)). \tag{3.6} \]

\begin{proof}

We proceed with a direct computation of the derivative:

\[
\frac{\partial}{\partial s} \phi_s (\exp(tH_0)) = \int_K \frac{\partial}{\partial s} \left[ e^{(i(s \alpha) - \rho)H(\exp(tH_0) k)} \right]  dk.
\]

The derivative of the integrand is:
\[
\frac{\partial}{\partial s} \left[ e^{(i(s \alpha) - \rho)H(\exp(tH_0) k)} \right] = i\alpha (H(\exp(tH_0) k)) \cdot e^{(i(s \alpha) - \rho)H(\exp(tH_0) k)}.
\]

Thus,
\[
\frac{\partial \phi_s}{\partial s} (\exp(tH_0)) = i \int_K \alpha (H(\exp(tH_0) k)) \cdot e^{(i(s \alpha) - \rho)H(\exp(tH_0) k)}  dk.
\]

\[
\frac{\partial^2}{\partial s^2} \phi_s(\exp(tH_0)) = -\int_K [\alpha (H(\exp(tH_0) k))]^2 \cdot e^{(i(s \alpha) - \rho)H(\exp(tH_0) k)} dk
\]

\[
\frac{\partial^4}{\partial s^4} \phi_s(\exp(tH_0)) = \int_K [\alpha (H(\exp(tH_0) k))]^4 \cdot e^{(i(s \alpha) - \rho)H(\exp(tH_0) k)}  dk
\]

Final inequality:

\[
\biggr| \frac{\partial^4}{\partial s^4} \phi_s(\exp(tH_0)) \biggr| = \biggr| \int_K [\alpha (H(\exp(tH_0) k))]^4 \cdot e^{(i(s \alpha) - \rho)H(\exp(tH_0) k)}   dk \biggr| \]
\[\leq
\int_K [\alpha (H(\exp(tH_0) k))]^4  \biggr|e^{- \rho H( \exp(tH_0) k)} \biggr| dk
\]
\[\leq \int_K (|H_\alpha| \cdot |H_0|)^2 t^2 [\alpha(H(\exp(tH_0) k))]^2 e^{- \rho H( \exp(tH_0) k)} dk = (|H_\alpha| \cdot |H_0|)^2 t^2 \biggr| \frac{\partial^2}{\partial s^2} \phi_{s}(t) |_{s=0}\biggr|
\]
The last inequality follows from Lemma 2.1: $ |\alpha(H(\exp(tH_0)k))| \leq (|H_\alpha| \cdot |H_0|) t.$

The second derivative is obtained by analogous computation of (3.4) combined with Lemma 2.1:
\[
-\frac{\partial^2}{\partial \lambda^2} \phi_{\lambda \alpha}(\exp(tH_0)) \biggr|_{\lambda=0} = \int_K \left[\alpha( H(\exp(tH_0) k))\right]^2 e^{-\rho H(\exp(tH_0) k)}  dk
\]
\[\leq \int_K  (t \| H_0 \| \| H_\alpha \| )^2 e^{-\rho H(\exp(tH_0) k)}  dk = (t \| H_0 \|  \| H_\alpha \|)^2 \varphi_{0} (\exp(tH_0)). \]

\end{proof}

\subsection{Non-Euclidean Mean, Variance, and Scaling}

We define, coherently with the Euclidean case, the mean and variance associated with a \(M\)-valued random variable \(Z\) defined on some probability space \((\Omega,\mathcal{F},\mathbb{P})\) satisfying assumption [R].

The \emph{analogue} of the characteristic function (of the second kind with the terminology of \cite{Karpelevich59}) writes:

\begin{equation}
\forall\lambda\in\mathbb{R},\ \Phi_{Z}(\lambda)=\frac{\hat{f}_{Z}(\lambda)}{\hat{f}_{Z}(0)}. \tag{3.7}
\end{equation}

Again the normalization guarantees that \(\Phi_{Z}(0)=1\) as in the Euclidean case. We also emphasize that, since we assumed the density to be radial, it follows from (3.1) and (3.4) that for all \(m=2j+1,j\in\mathbb{N}\),

\begin{equation}
\partial_{\lambda}^{m}\Phi_{Z}(\lambda)|_{\lambda=0}=0. \tag{3.8}
\end{equation}

In other terms, the \emph{odd} moments of the random variable are \(0\).

From the above definition we then define the \emph{analogue} of the variance as

\begin{equation}
V_{Z}=-\partial_{\lambda}^{2}\Phi_{Z}(\lambda)|_{\lambda=0}. \tag{3.9}
\end{equation}

\textbf{Proposition 3.3 (Variance of the sum).} Let \(Z_{1}\) and \(Z_{2}\) be two \(M\)-valued independent random variables with radial densities \(f_{Z_{1}},f_{Z_{2}}\in C_{0}(K \backslash G / K)\) w.r.t. the Riemannian volume {\color{brown}(2.1)} of \(M\). It then holds that

\[
V_{Z_{1}\oplus Z_{2}}=V_{Z_{1}}+V_{Z_{2}}.
\]

\begin{proof}
    The proof coincides with that in our previous paper \cite{KonMen25}.
\end{proof}

\textbf{Proposition 3.4 (scaled variables):}
\begin{align*}
\forall \eta \in (0, +\infty), \quad f_{\epsilon \otimes Z} (\eta) = \frac{1}{\epsilon}f_Z\biggr(\frac{\eta}{\epsilon}\biggr) \frac{(\sinh \frac{\eta}{\epsilon})^{m_\alpha} (\sinh (2\frac{\eta}{\epsilon}))^{m_{2\alpha}}}{(\sinh \eta)^{m_\alpha} (\sinh (2\eta))^{m_{2\alpha}}} \tag{3.10}
\end{align*}

\begin{proof}

$$\mu^{\epsilon}(B(r)) = \int_{G} I_{B(r)}\{\exp(\epsilon \log(g)\} f_Z(g) dg$$
In polar coordinates (with normalization constant $C$ for the probability measure): 
$$\mu^{\epsilon}(B(r)) = C\int_{\mathfrak{a}^+} I\{\epsilon a \leq r\} f_Z(a) J(a) da = C \int_{0}^{+\infty} I\{\epsilon \eta \leq r\} f_Z(\eta)(\sinh \eta)^{m_\alpha} (\sinh (2\eta))^{m_{2\alpha}} d\eta$$
Let's make a change of coordinates $\hat \eta = \epsilon \eta$, $d\eta = \frac{1}{\epsilon}d\hat\eta$:
$$= C \int_{0}^{+\infty} I\{\hat \eta \leq r\} \frac{1}{\epsilon}f_Z\biggr(\frac{\hat\eta}{\epsilon}\biggr)(\sinh \frac{\hat\eta}{\epsilon})^{m_\alpha} (\sinh (2\frac{\hat\eta}{\epsilon}))^{m_{2\alpha}} d\hat\eta$$
$$= C \int_{0}^{+\infty} I\{\hat \eta \leq r\} \frac{1}{\epsilon}f_Z\biggr(\frac{\hat\eta}{\epsilon}\biggr) \frac{(\sinh \frac{\hat\eta}{\epsilon})^{m_\alpha} (\sinh (2\frac{\hat\eta}{\epsilon}))^{m_{2\alpha}}}{(\sinh \hat\eta)^{m_\alpha} (\sinh (2\hat\eta))^{m_{2\alpha}}} (\sinh \hat\eta)^{m_\alpha} (\sinh (2\hat\eta))^{m_{2\alpha}} d\hat\eta$$
Thus, density of scaled variable has shape:
\begin{align*}
f_{\epsilon \otimes Z} (\eta) = \frac{1}{\epsilon}f_Z\biggr(\frac{\eta}{\epsilon}\biggr) \frac{(\sinh \frac{\eta}{\epsilon})^{m_\alpha} (\sinh (2\frac{\eta}{\epsilon}))^{m_{2\alpha}}}{(\sinh \eta)^{m_\alpha} (\sinh (2\eta))^{m_{2\alpha}}} 
\end{align*}

\end{proof}

\textbf{Proposition 3.5 (Variance for the random walks).}
    Let $Z$ satisfy [R]. Set for $\varepsilon > 0$, $Z_\varepsilon := \varepsilon \otimes Z$. It then holds that there exists $C > 0$ s.t.
    \begin{equation}
        \tag{3.11}
        V_{Z_{\varepsilon}} < C \varepsilon^2.
    \end{equation}

\begin{proof}
    Here we use notation $\varphi_{\lambda \alpha}(\eta) = \varphi_{\lambda \alpha}(\exp(\eta H_0)), \quad  f_Z(\eta) = f_Z(\exp(\eta H_0)), \quad f_{Z_{\varepsilon}}(\eta) = f_{Z_{\varepsilon}}(\exp(\eta H_0))$

    where $\alpha(H_0) = 1, \ \lambda \in \mathbb{R}$

    $$V_{Z_\epsilon} = -\frac{\partial^2}{\partial {\lambda}^2} C \int_0^{+\infty} f_{Z_\epsilon} (\eta) \varphi_{\lambda \alpha}(\eta) (\sinh \eta)^{m_\alpha} (\sinh (2\eta))^{m_{2\alpha}} d\eta \biggr|_{\lambda=0}$$

    $$= -\frac{\partial^2}{\partial {\lambda}^2} C \int_0^{+\infty} \frac{1}{\epsilon}f_Z\biggr(\frac{\eta}{\epsilon}\biggr) (\sinh \frac{\eta}{\epsilon})^{m_\alpha} (\sinh (2\frac{\eta}{\epsilon}))^{m_{2\alpha}}\varphi_{\lambda \alpha}(\eta) d\eta \biggr|_{\lambda=0}$$
    
    $$= C \int_0^{+\infty} f_Z(\widetilde\eta) (\sinh \widetilde\eta)^{m_\alpha} (\sinh (2\widetilde\eta))^{m_{2\alpha}} \left( - \frac{\partial^2}{\partial {\lambda}^2}\varphi_{\lambda \alpha}(\epsilon \widetilde\eta) \biggr|_{\lambda=0} \right) d\widetilde\eta$$

    $$\leq C \int_0^{+\infty} f_Z(\widetilde\eta) (\sinh \widetilde\eta)^{m_\alpha} (\sinh (2\widetilde\eta))^{m_{2\alpha}} (\epsilon \widetilde\eta)^2 \varphi_{0}(\epsilon \widetilde\eta) d\widetilde\eta$$
    
    $$\leq \epsilon^2 \widetilde C$$

    Applying (3.10) gives the second equality via cancellation. The first inequality follows from (3.5), and the second from the normalization 
    \[
    C \int_0^{+\infty} f_Z(\widetilde\eta) (\sinh \widetilde\eta)^{m_\alpha} (\sinh (2\widetilde\eta))^{m_{2\alpha}} (\epsilon \widetilde\eta)^2 d\widetilde\eta = 1
    \]
    combined with $\varphi_{0}(\epsilon \widetilde\eta) \leq 1$.
\end{proof}

\textbf{Proposition 3.6 (Variance for the random walks).}

In particular, choosing \(\varepsilon = \frac{1}{\sqrt{N}}\), the above control can be specified to derive that with the notation of (2.5):

\begin{equation}
\tag{3.12}
V_{S_N} = V_{\oplus_{j=1}^N \frac{1}{\sqrt{N}} \otimes Z^j} \xrightarrow[N]{} t := \frac{1}{n} \int_0^{+\infty} \eta^2 \mu_{Z,R}(d\eta),
\end{equation}

which is precisely the asymptotic variance appearing in Theorems~2.2 and~2.3. Furthermore, there exists \(C \geq 1\) s.t.

\begin{equation}
\tag{3.13}
|V_{S_N} - t| \leq \frac{C}{N},
\end{equation}




\begin{proof}

Let us choose $H', H_0 \in \mathfrak{a}$ such that: $\| H' \| = 1$, $\alpha(H_0) = \langle H_\alpha, H_0 \rangle = 1$, $\langle \| H_\alpha\| H', \| H_0\| H' \rangle = 1$, thus we get $\| H_\alpha\| \| H_0\| = 1$. For $\eta$ in bounded domain we can write the Taylor expansion with respect to $\epsilon \to 0$ (here we use Lemma 3.1):

$$ \varphi_{\lambda \alpha}(\exp(\epsilon \eta \| H_0 \| H')) = 1 + \frac{1}{2!} \left( -\frac{\| \lambda \alpha\|^2 + \| \rho \|^2}{n} \right) (\epsilon \| H_0 \| \eta)^2 + O_{\eta, \lambda}(\epsilon^4)
$$

Differentiate this expansion by $\lambda$:

$$ - \frac{\partial^2 }{\partial \lambda^2}\varphi_{\lambda \alpha}(\exp(\epsilon \eta \| H_0 \| H')) = -\frac{\partial^2 }{\partial \lambda^2} \left[1  + \frac{1}{2!} \left( -\frac{\lambda^2 \| \alpha\|^2 + \| \rho \|^2}{n} \right) (\epsilon \| H_0 \| \eta)^2 + O_{\eta, \lambda}(\epsilon^4) \right]
$$

\[ = -\frac{1}{2} \left( -\frac{2 \| \alpha \|^2}{n} \right) (\epsilon \| H_0 \|  \eta)^2 + O_{\eta, \lambda}(\epsilon^4) =  \frac{(\epsilon \eta)^2 }{n} + O_{\eta, \lambda}(\epsilon^4), \tag{3.14}
\]


$$V_{Z_{\epsilon}} (0) = \frac{ C \int_0^{+\infty} f_Z(\eta) (\sinh \eta)^{m_\alpha} (\sinh (2\eta))^{m_{2\alpha}} \left[-\partial_{\lambda}^2 \phi_{\lambda}(\epsilon \eta) \biggr|_{\lambda=0} \right] d\eta}
{C \int_0^{+\infty} f_Z(\eta) (\sinh \eta)^{m_\alpha} (\sinh (2\eta))^{m_{2\alpha}}\varphi_{0}(\epsilon \eta) d\eta}$$ 
$$= \frac{C \int_0^{+\infty} f_Z(\eta) (\sinh \eta)^{m_\alpha} (\sinh (2\eta))^{m_{2\alpha}} [\epsilon^2 \eta^2 \frac{1}{n} + O_{\eta}(\epsilon^4)] d\eta}
{C \int_0^{+\infty} f_Z(\eta) (\sinh \eta)^{m_\alpha} (\sinh (2\eta))^{m_{2\alpha}} d\eta + O(\epsilon^2)}$$
$$= \frac{\epsilon^2}{n} \frac{ C \int_0^{+\infty} f_Z(\eta) (\sinh \eta)^{m_\alpha} (\sinh (2\eta))^{m_{2\alpha}} \eta^2 d\eta}{C \int_0^{+\infty} f_Z(\eta) (\sinh \eta)^{m_\alpha} (\sinh (2\eta))^{m_{2\alpha}} d\eta} + O(\epsilon^4)$$
$$= \frac{\epsilon^2}{n} \int_0^{+\infty} \eta^2 \mu_{Z} (d\eta) + O(\epsilon^4)$$

Second equality we've got substituted (3.14) in formula.
\end{proof}

\section{Proof of the Theorem 2.3}



\[\hat{f}_{\frac{1}{\sqrt{N}}\otimes Z}(\lambda) = C \int_{0}^{+\infty} f(s) \phi_{-\lambda \alpha}(\frac{s}{\sqrt{N}}) \sinh(s)^{m_{\alpha}} \sinh(2s)^{m_{2\alpha}} ds\]

Denote $\hat{t} := \frac{1}{n \| \alpha\|^2}\int \eta^2 \mu_Z(d \eta) := \frac{t}{\| \alpha\|^2}$,
 $\rho := \frac{1}{2}m_{\alpha} + m_{2\alpha}$, $\hat{\rho} := (\frac{1}{2}m_{\alpha} + m_{2\alpha})\alpha$.
\begin{align*}
&f_{S_N}(\eta) - \Psi(\hat{t}, \eta) \nonumber\\
&= C \int_0^{+\infty} [\hat{f}_{S_N}(\lambda) - \hat{\Psi}(\hat{t}, \lambda)] 
   \varphi_\lambda \left( \eta \right) |c(\lambda)|^{-2}  d\lambda \nonumber \\
&= C \int_0^{+\infty} \left[ \prod_{j=1}^N \hat{f}_{1/\sqrt{N} \otimes Z^j}(\lambda) 
   - \exp \left( -\frac{(\| \hat{\rho} \|^2 + \| \lambda \alpha\|^2)\hat{t}}{2} \right) \right]
   \varphi_\lambda \left( \eta \right) |c(\lambda)|^{-2}  d\lambda \nonumber \\
&= C \int_0^{+\infty} \left[ \prod_{j=1}^N \hat{f}_{1/\sqrt{N} \otimes Z^j}(\lambda) 
   - \exp \left( -\frac{(\rho^2 + \lambda^2)t}{2} \right) \right]
   \varphi_\lambda \left( \eta \right) |c(\lambda)|^{-2}  d\lambda \nonumber \\
&= C \int_0^{+\infty} (I_{\lambda \leq D_N} + I_{\lambda > D_N}) 
   \left[ \left( \hat{f}_{1/\sqrt{N} \otimes Z}(\lambda) \right)^N 
   - \exp \left( -\frac{(\rho^2 + \lambda^2)t}{2} \right) \right] \nonumber \\
&\quad \times \varphi_\lambda \left( \eta \right) |c(\lambda)|^{-2}  d\lambda \nonumber
=: (\mathcal{B}_N + \mathcal{T}_N)(\eta). \tag{4.1}
\end{align*}



\textbf{Contribution of the Central Part.} First, we consider the difference between two Fourier transforms in the central part, i.e., for \(\lambda \leq D_{N}\). Specifically, we write (for details see \cite{KonMen25} page 435)

\begin{align*}
&\left| \left( \widehat{f}_{1/\sqrt{N} \otimes Z}(\lambda) \right)^{N} - \exp\left( - \frac{( \rho ^{2} + \lambda^{2})t}{2} \right) \right| \nonumber \\
&\leq \exp\left( -\frac{( \rho ^{2} + \lambda^{2})t}{2} \right) \left| N \ln\left( \widehat{f}_{1/\sqrt{N} \otimes Z}(\lambda) \right) + \frac{(\lambda^{2} +  \rho ^{2})t}{2} \right| \nonumber \\
&\qquad \times \exp\left( \left| N \ln\left( \widehat{f}_{1/\sqrt{N} \otimes Z}(\lambda) \right) + \frac{(\lambda^{2} +  \rho ^{2})t}{2} \right| \right),  \tag{4.2}
\end{align*}
\begin{align*}
&\left|N\ln(\widehat{f}_{1/\sqrt{N}\otimes Z}(\lambda))+\frac{(\lambda^2+  \rho ^2)t}{2}\right|\\
&\quad\leq\left|N\ln(\widehat{f}_{1/\sqrt{N}\otimes Z}(0))+\frac{ \rho ^2}{2}t\right|+ \frac{\lambda^2}{2}|t-V_{S_N}|\\
&\qquad\qquad+N\bigg(\sum_{r=2}^{\infty}r^{-1}\bigg|1-\frac{ \widehat{f}_{1/\sqrt{N}\otimes Z}(\lambda)}{\widehat{f}_{1/\sqrt{N}\otimes Z}(0)} \bigg|^r+\bigg|-1+\frac{\widehat{f}_{1/\sqrt{N}\otimes Z}(\lambda)}{ \widehat{f}_{1/\sqrt{N}\otimes Z}(0)}+\frac{\lambda^2}{2}V_{1/\sqrt{N}\otimes Z}\bigg| \bigg)\\
&=:\mathcal{R}_{1,N}+\mathcal{R}_{2,N}+\mathcal{R}_{3,N}.
\tag{4.3}
\end{align*}

We can write, expanding to the order four,

\begin{equation}
\frac{\widehat{f}_{1/\sqrt{N} \otimes Z}(\lambda)}{\widehat{f}_{1/\sqrt{N} \otimes Z}(0)} = 1 - \frac{\lambda^{2}}{2} V_{1/\sqrt{N} \otimes Z} + \frac{\lambda^{4}}{6} \int_{0}^{1} d\delta (1 - \delta)^{3} \frac{\widehat{f}^{(4)}_{1/\sqrt{N} \otimes Z}(\delta \lambda)}{\widehat{f}_{1/\sqrt{N} \otimes Z}(0)} \tag{4.4}.
\end{equation}





Let us remind (3.5):

\begin{align*}
\max_{\delta \in [0,1]} \biggl| \varphi^{(4)}_{\delta\lambda} ( \eta ) \biggr| \leq C \eta^2 \biggl| \varphi^{(2)}_0 ( \eta ) \biggr|
\end{align*}

and, therefore,

\begin{equation}
|R_N| \leq \mathbb{E} \left[ |Z_{1/\sqrt{N}, R}|^2 \frac{\varphi_0^{(2)} (Z_{1/\sqrt{N}, R})}{\hat{f}_{1/\sqrt{N} \otimes Z}(0)} \right],
\tag{4.5}
\end{equation}

where \( Z_{1/\sqrt{N}, R} \) is a random variable with the distribution law \( \mu_{1/\sqrt{N} \otimes Z, R}\). We have

\begin{align*}
&\mathbb{E}\left[|Z_{1/\sqrt{N},R}|^{2}\frac{|\varphi^{(2)}_{0}(Z_{1/ \sqrt{N},R})|}{\widehat{f}_{1/\sqrt{N}\otimes Z}(0)}\right] \\
&\quad=\Omega_{n-1}\int_{0}^{+\infty}\eta^{2}\frac{|\varphi^{(2)}_{0}(\eta)|}{\widehat{f}_{1/\sqrt{N}\otimes Z}(0)}f_{1/\sqrt{N}\otimes Z} \left(\eta\right)(\sinh \eta)^{m_\alpha} (\sinh (2\eta))^{m_{2\alpha}}\,d\eta \\
&\quad\underset{\mathclap{(3.18)}}{=}\Omega_{n-1}\int_{0}^{+\infty}\eta^{2}\frac{| \varphi^{(2)}_{0}(\eta)|}{\widehat{f}_{1/\sqrt{N}\otimes Z}(0)}\sqrt{N }\,f_{Z}\left(\sqrt{N}\,\eta\right) \\
&\quad\qquad\times \frac{(\sinh \sqrt{N} \eta)^{m_\alpha} (\sinh (\sqrt{N} 2\eta))^{m_{2\alpha}}}{(\sinh \eta)^{m_\alpha} (\sinh (2\eta))^{m_{2\alpha}}}(\sinh \eta)^{m_\alpha} (\sinh (2\eta))^{m_{2\alpha}}\,d\eta \\
&\quad=\Omega_{n-1}\frac{1}{\widehat{f}_{1/\sqrt{N}\otimes Z}(0)}\int_{0}^{+\infty}\frac{\widetilde{\eta}^{2}}{N}\Big{|}\varphi^{(2)}_{0}\left(\frac{\widetilde{\eta}}{N^{1/2}}\right)\Big{|} \\
&\quad\qquad\times f_{Z}\left(\widetilde{\eta}\right)(\sinh \widetilde\eta)^{m_\alpha} (\sinh (2\widetilde\eta))^{m_{2\alpha}}\,d\widetilde{\eta} \\
&\quad\leq C\Omega_{n-1}\frac{1}{\widehat{f}_{1/\sqrt{N}\otimes Z}(0)} \int_{0}^{+\infty}\frac{\widetilde{\eta}^{4}}{N^{2}}\Big{|}\varphi_{0}\left( \frac{\widetilde{\eta}}{N^{1/2}}\right)\Big{|} \\
&\quad\qquad\times f_{Z}\left(\widetilde{\eta}\right)(\sinh \widetilde\eta)^{m_\alpha} (\sinh (2\widetilde\eta))^{m_{2\alpha}}\,d\widetilde{\eta} \\
&\quad\leq\frac{C}{N^{2}} \tag{4.6}
\end{align*}


Therefore, there exists a \( C \) such that
\begin{equation} 
|R_N| \leq \frac{C}{N^2}. \tag{4.7}
\end{equation}

Proceeding as in \cite{KonMen25} (pages 436--437), we derive:
\begin{equation} 
\mathcal{R}_{2,N} \leq C \frac{\lambda^2}{N}. \tag{4.8}
\end{equation}
\begin{equation}
\mathcal{R}_{3,N} \leq C \lambda^4 N^{-1}. \tag{4.9}
\end{equation}

It remains to estimate the contribution of \( \mathcal{R}_{1,N} \) into (4.3). We have
\begin{align*}
\mathcal{R}_{1,N} &= \left| N \ln (\hat{f}_{1/\sqrt{N} \otimes Z}(0)) + \frac{\rho^2}{2}t \right| \\
&= \left| N \ln \left( 1 - \left( 1 - \hat{f}_{1/\sqrt{N} \otimes Z}(0) \right) \right) + \frac{\rho^2}{2}t \right| \\
&= \left| -N(1 - \hat{f}_{1/\sqrt{N} \otimes Z}(0)) + \frac{\rho^2}{2}t \right. \\
&\quad \left. - N(1 - \hat{f}_{1/\sqrt{N} \otimes Z}(0)) \sum_{r=2}^{\infty} r^{-1}(1 - \hat{f}_{1/\sqrt{N} \otimes Z}(0))^{r-1} \right|. \tag{4.10}
\end{align*}




Using definition of $\hat{f}_{1/\sqrt{N}\otimes Z}$: 

\[
\hat{f}_{1/\sqrt{N}\otimes Z}(0) = C \int_0^{+\infty} \varphi_0 \left(\eta \right) f_{\frac{1}{\sqrt{N}} \otimes Z} (\eta) (\sinh \eta)^{m_\alpha} (\sinh (2\eta))^{m_{2\alpha}}  d \eta \]

\[= C \int_0^{+\infty}  \varphi_0 \left(\eta \right)  \sqrt{N} f_Z\biggr(\sqrt{N} \eta\biggr) \frac{(\sinh \sqrt{N} \eta)^{m_\alpha} (\sinh (2 \sqrt{N}\eta))^{m_{2\alpha}}}{(\sinh \eta)^{m_\alpha} (\sinh (2\eta))^{m_{2\alpha}}} (\sinh \eta)^{m_\alpha} (\sinh (2\eta))^{m_{2\alpha}}  d \eta 
\]
\[= C \int_0^{+\infty}  \varphi_0 \left(\eta \right) f_Z\left(\sqrt{N}\eta \right) (\sinh \sqrt{N} \eta)^{m_\alpha} (\sinh (2 \sqrt{N}\eta))^{m_{2\alpha}} d \left(\sqrt{N}\eta \right)
\]
\[= C \int_0^{+\infty}  \varphi_0 \left( \frac{\widetilde\eta}{\sqrt{N}} \right) f_Z(\widetilde\eta) (\sinh \widetilde\eta)^{m_\alpha} (\sinh (2\widetilde \eta))^{m_{2\alpha}} d \widetilde \eta. \tag{4.11}
\]

Since $\mu_{\frac{1}{\sqrt{N}} \otimes Z}$ is a probability measure, the following equality holds:
\[C \int_0^{+\infty} f_{\frac{1}{\sqrt{N}} \otimes Z} (\eta) (\sinh \eta)^{m_\alpha} (\sinh (2\eta))^{m_{2\alpha}}  d \eta = 1. \tag{4.12} \]
From (4.11) and (4.12) it follows:
\[
\hat{f}_{1/\sqrt{N}\otimes Z}(0) - 1 = C \int_0^{+\infty} f_z \left( \tilde{\eta} \right) \left( \varphi_0 \left(  \frac{\tilde{\eta}}{\sqrt{N}} \right) - 1 \right) (\sinh \eta)^{m_\alpha} (\sinh (2\eta))^{m_{2\alpha}}   d\tilde{\eta}.
\]

Choosing now \( g(\zeta) = \varphi_0 (\zeta) \), we perform the Taylor expansion up to the fourth order, which gives

\[
\hat{f}_{1/\sqrt{N}\otimes Z}(0) - 1 = \Omega_{n-1} \int_0^{+\infty} f_z \left( \tilde{\eta} \right) 
 \left( g^{(1)}(0) \frac{\tilde{\| H_0 \| \eta}}{N^{1/2}} + \frac{1}{2} g^{(2)}(0) \frac{(\| H_0 \|\tilde{\eta})^2}{N} + \frac{1}{6} g^{(3)}(0) \frac{(\| H_0 \| \tilde{\eta})^3}{N^{3/2}} \right. \\ \]
 \[
 \left. + \frac{1}{6} \int_0^1 g^{(4)} (\gamma \tilde{\eta}) (1 - \gamma)^3  d\gamma \frac{(\| H_0 \| \tilde{\eta})^4}{N^2} \right) (\sinh \eta)^{m_\alpha} (\sinh (2\eta))^{m_{2\alpha}}  d\tilde{\eta}, \tag{4.13}
\]

where \( g^{(i)} \) stands for the \( i \)-th derivative of \( g \). Now it remains to calculate the derivatives of \( g \). Using g is an even function and Lemma 3.1, it gives:

\[ g^{(1)}(0) = 0, \ g^{(3)}(0) = 0, \ g^{(2)}(0) = -\frac{\rho^2 \| \alpha \|^2}{n}. \tag{4.14} \]

Substitute (4.14) in (4.13):

\[
\begin{aligned}
\hat{f}_{1/\sqrt{N}\otimes Z}(0) - 1 = C \int_0^{+\infty} f_z \left( \tilde{\eta} \right) 
& \left( - \frac{1}{2}\frac{\rho^2}{n}  \frac{\tilde{\eta}^2}{N} + \frac{1}{6} \int_0^1 g^{(4)} (\gamma \tilde{\eta}) (1 - \gamma)^3  d\gamma \frac{(\| H_0 \| \tilde{\eta})^4}{N^2} \right) \times\\ \\
&\times (\sinh \widetilde \eta)^{m_\alpha} (\sinh (2 \widetilde \eta))^{m_{2\alpha}}  d\tilde{\eta},
\end{aligned}
\]

Thus, we have established that  
\[
1 - \hat{f}_{1/\sqrt{N} \otimes Z}(0) = \frac{\rho^2}{n} \frac{1}{2N} \int_0^{+\infty} \tilde{\eta}^2 \mu_{Z,R}(d\widetilde\eta) + O\left(\frac{1}{N^2}\right),
\]  

which finally gives in (4.10) the expression

\[
\mathcal{R}_{1,N} = -\frac{\rho^2}{2} \frac{1}{n} \int_0^{+\infty} \tilde{\eta}^2  \mu_{Z,R}(d\eta) + O\left(\frac{1}{N}\right) + \frac{\rho^2}{2}t.
\]

Ultimately, this gives from the very definition of \( t \) (see, e.g., (3.13)) that

\begin{equation} \tag{4.15}
\mathcal{R}_{1,N} = O\left(\frac{1}{N}\right).
\end{equation}

Therefore, plugging (4.15), (4.9), (4.8) into (4.3) and (4.2) we derive, choosing \( D_N = N^{\frac{1}{4}} \) so that in particular \(\lambda^4 / N \leq 1\):


\[|B_N(\eta)| \leq C \int_{0 < \lambda \leq D_N} \exp(-\frac{(\rho^2 + \lambda^2)t}{2}) \frac{1}{N} \left( 1 + \lambda^4 \right) \exp\left( \frac{1}{N} \left( 1 + \lambda^4 \right) \right) |\varphi_\lambda (\eta)||c(\lambda)|^{-2} d\lambda\]

\[\leq \frac{C}{t^2 N} \int_{0 < \lambda \leq D_N} \exp(-\frac{(\rho^2 + \lambda^2)t}{4}) |\varphi_\lambda (\eta)||c(\lambda)|^{-2} d\lambda\]

\[\leq \frac{C}{t^2 N} \int_{\mathbb{R}^+} \exp(-\frac{(\rho^2 + \lambda^2)t}{4}) |\varphi_\lambda (\eta)||c(\lambda)|^{-2} d\lambda.\]

We can estimate this integral using following inequalities from \cite{Platonov99} and \cite{helgason1984groups}:
\[
|\varphi_\lambda (\eta)| \leq  1, \quad |c(\lambda)|^{-2} \leq C (1 + \lambda^{n}), \tag{4.16}
\]

from which we eventually derive:

\[
|B_N(\eta)| \leq \frac{C}{t^2 N} \left( \frac{1}{t^{\frac{1}{2}}} \wedge \frac{1}{t^{\frac{n}{2}}} \right). \tag{4.17}
\]

\textbf{Contribution of the tails, general case.} For $\lambda > N^{1/4}$, we write:

\begin{align}
\mathcal{T}_N(\eta) &\leqslant \left| \int_{N^{1/4}}^{+\infty} \left( \prod_{j=1}^N \hat{f}_{\frac{1}{\sqrt{N}}\otimes Z^j}(\lambda) \right) \varphi_\lambda(\eta)) |\mathbf{c}(\lambda)|^{-2} d\lambda \right| + \left| \int_{N^{1/4}}^{+\infty} \exp(-\frac{(\rho^2 + \lambda^2)t}{2}) \varphi_\lambda(\eta)) |\mathbf{c}(\lambda)|^{-2} d\lambda \right| \notag \\
&=: (\mathcal{T}_N^1 + \mathcal{T}_N^2)(\eta). \tag{4.18}
\end{align}

Let us first consider the term $\mathcal{T}_N^2$ which can be handled globally. We get again from (4.16) (similarly to (4.17)):

\begin{align}
\mathcal{T}_N^2(\eta) &\leqslant \exp(-\frac{N^{\frac{1}{2}}t}{4}) \int_0^{+\infty} \exp(-\frac{(\rho^2 + \lambda^2)t}{4}) |\varphi_\lambda(\eta))| |\mathbf{c}(\lambda)|^{-2} d\lambda \notag \\
&\leqslant \frac{C}{t^2 N} \left( \frac{1}{t^{\frac{1}{2}}} \wedge \frac{1}{t^{\frac{n}{2}}} \right), \tag{4.19}
\end{align}

which gives an upper bound similar to the one obtained for the bulk in (4.17).

Let us now turn to $\mathcal{T}_N^1(\eta)$, and split as follows:
\begin{align}
\mathcal{T}_N^1(\eta) &\leqslant \left| \int_{N^{\frac{1}{4}}}^{c_0 N^{\frac{1}{2}}} \left( \hat{f}_{\frac{1}{\sqrt{N}}\otimes Z}(\lambda) \right)^N \varphi_\lambda(\eta)) |\mathbf{c}(\lambda)|^{-2} d\lambda \right| \notag \\
&\quad + \left| \int_{c_0 N^{\frac{1}{2}}}^{\infty} \left( \hat{f}_{\frac{1}{\sqrt{N}}\otimes Z}(\lambda) \right)^N \varphi_\lambda(\eta)) |\mathbf{c}(\lambda)|^{-2} d\lambda \right| =: (\mathcal{T}_N^{11} + \mathcal{T}_N^{12})(\eta), \tag{4.20}
\end{align}

for some small enough constant $c_0$ to be specified. 
Now we want to estimate following term:

\begin{equation} \tag{4.21}
\hat{f}_{\frac{1}{\sqrt{N}}\otimes Z}(\lambda) = \int_0^\infty \varphi_\lambda \left( \frac{\eta}{\sqrt{N}} \right) \mu_{Z,R}(d\eta).
\end{equation}

We will use the result from \cite{helgason1984groups} page 252, where $\varphi_\lambda$ is represented by hypergeometric function: 

\[
a = \frac{1}{2} \left( \frac{1}{2} m_\alpha + m_{2\alpha} + \langle i \lambda, \alpha_0 \rangle \right),
\]

\[
b = \frac{1}{2} \left( \frac{1}{2} m_\alpha + m_{2\alpha} - \langle i \lambda, \alpha_0 \rangle \right),
\]

\[
c = \frac{1}{2} \left( m_\alpha + m_{2\alpha} + 1 \right).
\]
where, with \( \alpha_0 = \alpha / \langle \alpha, \alpha \rangle \),

Thus \( \varphi_\lambda \) is given by the hypergeometric function

\[
\varphi_\lambda (h) = F(a, b; c; z), \quad z = -\sinh^2 (\alpha (\log h)).
\]

Let's use three conventions to compute this function in coordinates:

$$\alpha(H_0) = 1, \ \lambda = s \alpha, \ h = \exp(t H_0), \ \eta' = \frac{\eta}{\sqrt{N}}$$

Thus:

\[
a = \frac{1}{2} \left( \frac{1}{2} m_\alpha + m_{2\alpha} + \langle i s \alpha, \alpha \rangle / \langle \alpha, \alpha \rangle \right) = 
\frac{1}{2} \left( \frac{1}{2} m_\alpha + m_{2\alpha} + i s \right),
\]

\[
b = \frac{1}{2} \left( \frac{1}{2} m_\alpha + m_{2\alpha} - \langle i s \alpha, \alpha \rangle / \langle \alpha, \alpha \rangle \right) = \frac{1}{2} \left( \frac{1}{2} m_\alpha + m_{2\alpha} - i s \right),
\]

\[
c = \frac{1}{2} \left( m_\alpha + m_{2\alpha} + 1 \right).
\]

$$
z = -\sinh^2 (\alpha (\log \exp(t H_0))) =  -\sinh^2 (t).
$$

\begin{equation}
    {}_2F_1(a, b; c; z) = \sum_{q=0}^{+\infty} \frac{(a)_q (b)_q}{(c)_q} \frac{z^q}{q!}, \tag{4.22}
\end{equation}

denoting for \( d \in \{a, b, c\} \),

\begin{equation}
    (d)_q = 
    \begin{cases} 
        1, & q = 0, \\
        \prod_{j=0}^{q-1} (d+j), & q \neq 0.
    \end{cases}
    \tag{4.23}
\end{equation}


\begin{align*}
\varphi_{\lambda \alpha}\left(\exp(\eta' H_0)\right) 
&= 1 + \sum_{q=1}^{\infty} \frac{\prod_{j=0}^{q-1} \left((\rho/2 + i\frac{\lambda}{2} + j)(\rho/2 - i\frac{\lambda}{2} + j)\right) \, (-1)^q \left(\sinh\left(\eta' \right)\right)^{2q}}{\prod_{j=0}^{q-1} \left(\frac{n}{2} + j\right) q!} \\
&= 1 + \sum_{q=1}^{\infty} \frac{\prod_{j=0}^{q-1} \left((\rho/2 + j)^2 + \frac{\lambda^2}{4}\right) \, (-1)^q \left(\sinh\left( \eta' \right)\right)^{2q}}{\prod_{j=1}^{q} (m + j) q!}, \quad m =  \frac{n}{2} - 1, \\
&= 1 + \sum_{q=1}^{\infty} \frac{\prod_{j=0}^{q-1} \left(\left(2\frac{\rho/2+j}{\lambda}\right)^2 + 1\right) \, (-1)^q \left(\lambda \sinh\left(\eta' \right)\right)^{2q}}{4^q \prod_{j=1}^{q} (m + j) q!} \\
&= 1 + \sum_{q=1}^{\infty} \frac{(-1)^q \left(\lambda \sinh\left(\eta'\right)\right)^{2q} \left[ 1 + \left( \frac{\rho}{\lambda} \right)^2 \right] \times \cdots \times \left[ 1 + \left( \frac{\rho+2q-2}{\lambda} \right)^2 \right] }{4^q q! (m + 1) (m + 2) \cdots (m + q)} 
\end{align*}
$$= 1 - \left(\frac{\lambda \sinh\left( \eta' \right)}{2} \right)^{2} \frac{1 + \left( \frac{\rho}{\lambda} \right)^2}{m + 1} + \left(\frac{\lambda \sinh\left(\eta'\right)}{2} \right)^{2} \frac{1 + \left( \frac{\rho}{\lambda} \right)^2}{m + 1} \times$$
$$
\times \left[ \sum_{q=2}^{\infty} \frac{(-1)^q \left(\lambda \sinh\left(\eta' \right)\right)^{2(q-1)} \left[ 1 + \left( \frac{\rho + 2}{\lambda} \right)^2 \right] \times \cdots \times \left[ 1 + \left( \frac{\rho+2q-2}{\lambda} \right)^2 \right] }{4^{q-1} q! (m + 2) \cdots (m + q)} \right],$$

\[
\frac{\left[1+\left(\frac{\rho + 2}{\lambda}\right)^{2}\right]\times\cdots\times \left[1+\left(\frac{\rho+2q-2}{\lambda}\right)^{2}\right]}{q!(m+2)\cdots(m+ q)} = \prod_{k=2}^{q} \frac{1 + (\frac{\rho + 2(k - 1)}{\lambda})^2}{k(m + k)}.
\]

Let's observe k-th term independently:

$$\frac{1 + (\frac{\rho + 2(k - 1)}{\lambda})^2}{k(m + k)} = \frac{\frac{1}{k^2} + \frac{1}{\lambda^2}(\frac{\rho^2}{k^2} + 4 \rho \frac{k - 1}{k} + (\frac{k - 1}{k})^2)}{1 + \frac{m}{k}} \leq 1,$$

\[
\prod_{k=2}^{q} \frac{1 + (\frac{\rho + 2(k - 1)}{\lambda})^2}{k(m + k)} \leq 1.
\]

Substitute inequality in formula:
$$
 \sum_{q=2}^{\infty} \frac{(-1)^q \left(\lambda \sinh\left(\eta\right)\right)^{2(q-1)} \left[ 1 + \left( \frac{\rho + 2}{\lambda} \right)^2 \right] \times \cdots \times \left[ 1 + \left( \frac{\rho+2q-2}{\lambda} \right)^2 \right] }{4^{q-1} q! (m + 2) \cdots (m + q)} 
$$

$$
 \leq \sum_{q=2}^{\infty} \frac{\left(\lambda \sinh\left(\eta' \right)\right)^{2(q-1)}}{2^{2(q-1)}}
 = \sum_{q=1}^{\infty} (\frac{\left(\lambda \sinh\left(\eta' \right)\right)}{2})^{2q} 
 = \sum_{q=1}^{\infty} (\frac{\left(\lambda \sinh\left(\eta' \right)\right)}{2^{1- \delta}})^{2q} (\frac{1}{2^\delta})^{2q}
$$

$$\leq \sum_{q=1}^{\infty}(\frac{1}{2^\delta})^{2q} =: q < 1,
$$

To get last inequality we use $\frac{\eta}{\sqrt{N}} \leq \frac{\bar{R}}{\sqrt{N}}$ because \(\mu_Z\) in (4.21) is compactly supported on \([0, \bar{R}]\) and decomposition of $\sinh\left(\frac{\eta}{\sqrt{N}}\right)$ in the heighborhood of zero:

$$
\frac{\lambda \sinh\left( \frac{\eta}{\sqrt{N}} \right)}{2^{1 - \delta}} 
\leq \frac{\lambda \cdot \sinh\left( \frac{\eta}{\sqrt{N}} \right)}{2^{1 - \delta}} 
\leq \frac{c_0 \sqrt{N} \cdot \eta [1+\delta]}{2^{1 - \delta} \sqrt{N}} 
\leq \frac{c_0 \bar{R} (1+\delta)}{2^{1 - \delta}}
\leq 1
$$
We need choose $c_0$ to hold the last equality (concrete values of $c_0$ see in (4.26)).

Finally, we get:

$$\varphi_{\lambda}\left(\exp \left(\frac{\eta}{\sqrt{N}} H_0 \right)\right) \leq 1 - A(1- q), \quad 
A = \left(\frac{\lambda \sinh\left(\eta' \right)}{2} \right)^{2} \frac{1 + \left( \frac{\rho}{\lambda} \right)^2}{m + 1}.$$

We then derive that on the considered range:

\begin{equation}
\varphi_\lambda \left(\frac{\eta}{\sqrt{N}}\right) \leq 1 - c\lambda^2 \left( \frac{\eta}{\sqrt{N}} \right)^2, \tag{4.24}
\end{equation}

where $c = \frac{1}{4(m + 1)} < \frac{1 + \left( \frac{\rho}{\lambda} \right)^2}{4(m + 1)}.$
Hence, from (4.21) and (4.24) we obtain:

\[
|\hat{f}_{\frac{1}{\sqrt{N}} \otimes Z}(\lambda)| \leq \left(1 - c \frac{\lambda^2}{N} \int_{0}^{\infty} \eta^2 \mu_{Z,R}(d\eta)\right), \quad |\hat{f}_{\frac{1}{\sqrt{N}} \otimes Z}(\lambda)|^N \leq \exp\left(-c\lambda^2 \int_{0}^{\infty} \eta^2 \mu_{Z,R}(d\eta)\right) = \exp(-c\lambda^2 t),
\]

up to a modification of \( c \) for this very last inequality. We then eventually get from (4.20), analogously to (4.19), that

\[
T_N^{11}(\eta) \leq C \exp\left(-c\frac{(-\rho^2 + N^{\frac{1}{2}})t}{2}\right) \int_{N^{\frac{1}{4}}}^{c_0 N^{\frac{1}{2}}} \exp\left(-c\frac{(\rho^2 + \lambda^2)t}{2}\right) |\varphi_\lambda (\eta)| |c(\lambda)|^{-2}  d\lambda
\]

\[
\leq \frac{C}{t^2 N} \left( \frac{1}{t^{\frac{1}{2}}} \wedge \frac{1}{t^{\frac{\eta}{2}}} \right). \tag{4.25}
\]

\textbf{It now remains to consider the case} $\lambda > c_0 N^{1/2}$ to handle $T_N^{12}(n)$ in (4.20).

\[
T_{N}^{12}(\eta) = \left| \int_{c_0 N^{\frac{1}{2}}}^{\infty} \left( \prod_{j=1}^{N} \hat{f}_{\frac{1}{\sqrt{N}} \otimes Z}(\lambda) \right) \varphi_{\lambda} \left( \eta \right) |c(\lambda)|^{-2} d\lambda \right| \leq \int_{c_0 N^{\frac{1}{2}}}^{\infty} \left| \hat{f}_{\frac{1}{\sqrt{N}} \otimes Z}(\lambda) \right|^{N} \left|\varphi_{\lambda} \left(\eta\right)\right| |c(\lambda)|^{-2} d\lambda 
\tag{4.26}
\]

Using equation (4.21),
\[
\left|\hat{f}_{\frac{1}{\sqrt{N}}\otimes Z}(\lambda) \right| \leq E \left[ \left|\varphi_{\lambda} \left( \frac{\eta}{\sqrt{N}} \right)\right| \right] = \int_{0}^{\eta_0} \varphi_{-\lambda} \left( \frac{\eta}{\sqrt{N}} \right) d \mu_Z(\eta) + \int_{\eta_0}^{\infty} \varphi_{-\lambda} \left( \frac{\eta}{\sqrt{N}} \right) d \mu_Z(\eta).
\]

Using the bound $|\varphi_{\lambda}(t)| \leq 1$, we obtain:

\[
\left| \hat{f}_{\frac{1}{\sqrt{N}}\otimes Z}(\lambda) \right| \leq \int_{0}^{\eta_0}  d\mu_Z(\eta') + \sup_{\eta' \geq \eta_0} \left| \varphi_{-\lambda} \left( \frac{\eta'}{\sqrt{N}} \right) \right| \int_{\eta_0}^{\infty} d\mu_Z(\eta') = \mu_Z(\eta' < \eta_0) + \sup_{\eta' \geq \eta_0} \left| \varphi_{-\lambda} \left( \frac{\eta'}{\sqrt{N}} \right) \right| p_0.
\]

This simplifies to ($p_0 = \mu_Z(\eta > \eta_0)$):

\[
\left| \hat{f}_{\frac{1}{\sqrt{N}}\otimes Z}(\lambda) \right| \leq (1 - p_0) + p_0 \cdot \sup_{\eta' \geq \eta_0} \left| \varphi_{-\lambda} \left( \frac{\eta'}{\sqrt{N}} \right) \right|.
\]

We now bound the supremum. For $\lambda \geq c_0 N^{1/2}$ and $\eta' \geq \eta_0$, the argument $t = \eta' / \sqrt{N}$ satisfies:
\[
\lambda t = \lambda \cdot \frac{\eta'}{\sqrt{N}} \geq (c_0 N^{1/2}) \cdot \frac{\eta_0}{\sqrt{N}} = c_0 \eta_0.
\]

We begin by recalling a result regarding the bound of the spherical function.

[Lemma 3.3 in \cite{Platonov99}]
For some constant $c > 0$, depending only on $M$, the inequality
\begin{equation} \label{eq:platonov_bound}
    1 - \varphi_\lambda(\exp(tH_0)) \geq c \tag{4.27}
\end{equation}
holds for all $\lambda t \geq 1$.



To apply this lemma, we select constants $c_0$ and $\eta_0$ such that $c_0 \eta_0 \geq 1$. Specifically, we set:
\[
c_0 = \frac{3}{2\bar{R}}, \quad \eta_0 = \frac{2\bar{R}}{3}. \tag{4.28}
\]

\[
\left| \varphi_{\lambda} \left( \frac{\eta'}{\sqrt{N}} \right) \right| \leq 1 - c =: \gamma_0.
\]

\[
\left| \hat{f}_{\frac{1}{\sqrt{N}}\otimes Z}(\lambda) \right| \leq (1 - p_0) + p_0 \gamma_0 =: \delta < 1
\]

Substitute in formula (4.25):

\[
T_{N}^{12}(\eta) \leq \delta^{N-2} \int_{c_0 N^{\frac{1}{2}}}^{\infty} \left| f_{\frac{1}{\sqrt{N}} \otimes Z}(\lambda) \right|^2 |\varphi_{\lambda}(\eta)| |c(\lambda)|^{-2} d\lambda
\]

We can now use the Plancherel equality, see \cite{helgason1984groups} page 454:

\begin{align*}
\mathcal{T}_{N}^{12}(\eta) &\leqslant C\delta^{N-2}\int_{0}^{\bar{R}}f_{\frac{1}{\sqrt{N}} \otimes_{Z}}^{2} \left(\eta\right)(\sinh \eta)^{m_\alpha} (\sinh (2\eta))^{m_{2\alpha}} d\eta \\
&\leqslant C\delta^{N-2}\int_{0}^{\frac{\bar{R}}{N^{\frac{1}{2}}}}Nf_{Z}^{2}  \left(
\frac{\sinh (N^{\frac{1}{2}}\eta)^{m_\alpha} (\sinh (2N^{\frac{1}{2}} \eta))^{m_{2\alpha}}}{(\sinh \eta)^{m_\alpha} (\sinh (2\eta))^{m_{2\alpha}}} 
\right)^2 (\sinh \eta)^{m_\alpha} (\sinh (2\eta))^{m_{2\alpha}} d\eta \\
&\leqslant C\delta^{N-2}N^{\frac{1}{2}}\int_{0}^{\bar{R}}f_{Z}^{2} \left(\eta\right) 
\left((\sinh \eta)^{m_\alpha} (\sinh (2\eta))^{m_{2\alpha}} \right)^2
 \sinh \left(\frac{\eta}{N^{\frac{1}{2}}}\right)^{-m_\alpha} \sinh \left(\frac{2\eta}{N^{\frac{1}{2}}}\right)^{-m_{2\alpha}}d\eta
\end{align*}

using as well (3.10) with $\varepsilon = N^{-\frac{1}{2}}$ for the second inequality. Hence, using $m_{\alpha} + m_{2 \alpha} = n - 1$

\[
\mathcal{T}_{N}^{12}(\eta) \leq C\delta^{N-2}N^{\frac{n}{2}}\int_{0}^{\bar{R}}f_{Z}^{2} \left(\eta\right) 
\left((\sinh \eta)^{m_\alpha} (\sinh (2\eta))^{m_{2\alpha}} \right)^2
\eta^{- (m_{\alpha} + m_{2\alpha}) } d\eta \leq C\delta^{N-2}N^{\frac{n}{2}}
\]

recalling that $\delta<1$ for the last inequality. We have thus established from the above control and (4.25) that,

\[
\mathcal{T}_{N}^{1}(\eta)\leqslant\frac{C}{N}\Big{(}\frac{1}{t^{2}}\Big{(} \frac{1}{t^{\frac{1}{2}}}\wedge\frac{1}{t^{\frac{n}{2}}}\Big{)}+1\Big{)}.
\]

From the above control and (4.19), (4.18) we thus derive:

\[
\mathcal{T}_{N}(\eta)\leqslant\frac{C}{N}\Big{(}\frac{1}{t^{2}}\Big{(}\frac{1}{ t^{\frac{1}{2}}}\wedge\frac{1}{t^{\frac{n}{2}}}\Big{)}+1\Big{)},
\]

which together with (4.17) and (4.1) completes the proof of the local limit Theorem 2.3.

\printbibliography

\end{document}